\documentclass[12pt]{amsart}

\usepackage[letterpaper, margin=1in]{geometry}

\usepackage[english]{babel}
\usepackage{amsmath}  
\usepackage{amssymb}  
\usepackage{amsfonts}
\usepackage{amsthm}     
\usepackage{physics}
\usepackage{graphicx, enumerate}  

\usepackage[colorlinks,linkcolor=blue,citecolor=blue,urlcolor=blue]{hyperref} 

\usepackage[backend=bibtex,isbn=false,url=false,doi=false]{biblatex} 

\newtheorem{theorem}{Theorem}[section] 
\newtheorem{lemma}[theorem]{Lemma}
\newtheorem*{lem*}{Lemma}

\newtheorem{proposition}[theorem]{Proposition}

\theoremstyle{definition}                  
\newtheorem{definition}[theorem]{Definition}
\newtheorem{conjecture}[theorem]{Conjecture}

\theoremstyle{remark}                       
\newtheorem{example}[theorem]{Example}  
\newtheorem*{exmp*}{Example}

\numberwithin{equation}{section}

\newcommand{\N}{\mathbb{N}}
\newcommand{\Z}{\mathbb{Z}}

\title[Random Walks on the Generalized Symmetric Group]{Random Walks on the Generalized Symmetric Group: Cutoff for the One-sided Transposition Shuffle}
\author{Yongtao Deng}
\author{Shi Jie Samuel Tan}
\date{}
\thanks{This work was partially supported by NSF Grant DMS \#2202017.}
\begin{document}

\maketitle

\begin{abstract}
    In this paper, we present a detailed proof for the exhibition of a cutoff for the one-sided transposition (OST) shuffle on the generalized symmetric group $G_{m,n}$. 
    Our work shows that based on techniques for $m \leq 2$ proven in \cite{M20}, we can prove the cutoff in total variation distance and separation distance for an unbiased OST shuffle on $G_{m,n}$ for any fixed $m \geq 1$ in time $n \log(n)$. We also prove the branching rules for the simple modules of $G_{m,n}$ and lay down some of the mathematical foundation for proving the conjecture for the cutoff in total variation distance for any general biased OST shuffle on $G_{m,n}$. 
\end{abstract}

\section{Introduction}

Random walks on finite groups have been intensively studied since the 1980s because they are good analogies to real-world problems.
Aldous and Diaconis gave a detailed introduction to the relationship between random walks on finite groups and the fields of representation theory and combinatorics \cite{MR0841111, MR0876954}. In their exposition, they also introduced the concept of mixing time and cutoff time for random walks.

In 2019, Bate, Connor, and Matheau-Raven introduced a new type of card shuffling called the \emph{one-sided transposition} (OST) shuffle which defines a particular random walk on the symmetric group.
They provided a formula for all the eigenvalues of this shuffle and proved a total variation cutoff for this shuffle at time $n \log(n)$ \cite{MR4312845}.

In 2020, Matheau-Raven's unpublished Ph.D. dissertation \cite{M20} presented random walks defined by the OST shuffle on the symmetric groups and the hyperoctahedral groups. 
The results included a cutoff time for the OST shuffle on these two groups. 
In the dissertation, Matheau-Raven shared open questions for random walks defined by the OST shuffle on the generalized symmetric groups $G_{m,n}$, and provided outlines of proofs for these questions.
In this paper, we extend some of the techniques from \cite{M20} to provide a formal proof that the cutoff time for the unbiased OST shuffle on the generalized symmetric group is $n \log(n)$. We also proved a conjecture from \cite{M20} regarding the branching rules for the simple modules of the generalized symmetric group and constructed operators that may prove to be useful for finding the cutoff time for the biased OST shuffle and the random transposition shuffle. 


\subsection{Random Walks on Finite Groups and Their Applications}

In 1988, Diaconis gave a detailed introduction on random walks on finite groups \cite{MR0964069}, and started a field on studying random walks defined by different card shuffling techniques. The early literature mostly focused on whether a deck can be well mixed with a type of shuffling technique, and how many shuffles are required to reach a well-mixed deck.

Aldous, Diaconis, and Shahshahani then discovered the cutoff phenomenon; in other words, they noticed that the probability distribution for the random walk reached a stationary distribution abruptly \cite{MR0770418, MR0626813}.
The cutoff phenomenon was later formalized by Aldous and Diaconis in \cite{MR0876954, MR1374011}. 
Although the cutoff phenomenon is found in limited examples, it is believed that this is widespread in random walks \cite{MR2023654}.

Because of the deep connection between random walks on finite groups and card shuffling, random walks have wide applications in gambling industries. 
On top of that, Diaconis also elucidated the connection between random walks and statistics \cite{MR0964069}.
Sinclair and Jerrum also related random walks on finite groups to Monte-Carlo Markov chains and hence to theoretical computer science problems \cite{MR1678570, sinclair2012algorithms}.

\subsection{Main Result}
In this work, we build on the theoretical foundation that Matheau-Raven has established on symmetric groups and the hyperoctahedral group, and we prove a theorem originally conjectured in Matheau-Raven's unpublished dissertation \cite{M20}. We first state the theorem informally in this section before stating it formally as Theorem \ref{thm:main_formal} in Section \ref{section:results}. All relevant background definitions and notations can be found in Section \ref{section:prelim}. 

\begin{theorem}[Conjectured in Section 4.5 in \cite{M20}]\label{thm:main}
The unbiased one-sided transposition shuffle on the generalized symmetric group $G_{m,n}$ approaches the uniform distribution after $n\log(n)$ shuffles.
\end{theorem}

Theorem \ref{thm:main} appears as Theorem 4.5.9 in the unpublished dissertation of Matheau-Raven, where an sketch for the proof was provided. In this paper, we formally prove the lower bound for the cutoff for the OST shuffle on $G_{m,n}$ by constructing a homomorphism from $G_{m,n}$ to $S_n$. We then connect the coupon collector's problem to the OST on $G_{m,n}$ to find the upper bound, allowing us to formally prove Theorem \ref{thm:main_formal} using both the lower and upper bounds (Propositions \ref{lem:lowerbound_OST_gensymgroup} and \ref{lem:upperbound_OST_gensymgroup}).

In addition, we proved Conjecture 4.5.4 in \cite{M20} that is informally stated as Theorem~\ref{conjecture:informal_branchingrules} in this section. We formally state it as Theorem~\ref{thm:branching_rules_gensymgroup} in Section~\ref{sec:branchingrules} and utilize the Littewood-Richardson Rule for $G_{m,n}$ stated in \cite{puttaswamaiah1969unitary} to prove it. 

\begin{theorem}[{\cite[Conjecture 4.5.4]{M20}}]\label{conjecture:informal_branchingrules}
    The branching rules for the simple modules of the generalized symmetric group $G_{m,n}$ are as follows:
    \begin{enumerate}
        \item The restricted representation of the simple module of $G_{m,n}$ that corresponds to the a particular partition $\lambda$ is the direct sum of simple modules of $G_{m,n-1}$ that corresponds to partitions that are contained within $\lambda$.
        \item The induced representation of the simple module of $G_{m,n}$ that corresponds to the a particular partition $\lambda$ is the direct sum of simple modules of $G_{m,n+1}$ that corresponds to partitions that contain $\lambda$.
    \end{enumerate}
\end{theorem}

On top of that, we proposed operators with certain properties  in Section~\ref{sec:future-work}. We believe that these operators and the properties that they possess, combined with Theorem~\ref{thm:branching_rules_gensymgroup}, can be used to tackle the following conjectures in \cite{M20} that are stated below and further elaborated in Section~\ref{sec:future-work}:

\begin{conjecture}[{\cite[Conjecture 4.5.5]{M20}}]\label{conjecture:lifting-eigenvectors_biased_OST}
    The eigenvalues for the biased one-sided transposition shuffle on the
generalized symmetric group $G_{m,n}$ are labelled by standard Young tableaux of certain partitions, and may be described by the technique of lifting eigenvectors.
\end{conjecture}

\begin{conjecture}[{\cite[Conjecture 4.5.6]{M20}}]\label{conjecture:lifting-eigenvectors_random}
    The eigenvalues for the random transposition shuffle on the
generalized symmetric group $G_{m,n}$ are labelled by standard Young tableaux of certain partitions, and may be described by the technique of lifting eigenvectors.
\end{conjecture}

\begin{conjecture}[{\cite[Conjecture 4.5.8]{M20}}]\label{conjecture:biasedOST}
    The biased one-sided transposition shuffle on the generalized symmetric group $G_{m,n}$ with weight function $w$ such that $w(j) = j^\alpha$ exhibits a cutoff in total variation distance at time $t_{n,\alpha} \log n$.
\end{conjecture}

\begin{conjecture}[{\cite[Conjecture 4.5.7]{M20}}]\label{conjecture:random}
    The random transposition shuffle on the generalized symmetric group $G_{m,n}$ exhibits a cutoff in total variation distance at time $t_{n,\alpha} \log n$.
\end{conjecture}

\subsection{Discussion of Related and Further work}
We now discuss some of the work done by others which is related to the one-sided transposition shuffle and random walks on the generalized symmetric group.

Ever since the one-sided transposition shuffle was proposed, Bate, Connor, and Matheau-Raven were able to prove that the one-sided transposition shuffle on the symmetric group exhibits a cutoff at $n \log n$ \cite{MR4312845}. Since then, Grinberg and Lafrenière have computed and analyzed the eigenvalues of the shuffling operators including the one-sided transposition shuffling operators on the symmetric group \cite{grinberg2023onesided}. They described a strong stationary time for the random-to-below shuffle, a generalized version of the one-sided transposition shuffle. Nestoridi and Peng introduced the one-sided $k$-transposition shuffle, a different generalization of the one-sided transposition shuffle, and studied the mixing time of this new shuffling technique on $S_n$ \cite{nestoridi2021mixing}. The one-sided $k$-transposition shuffle differs from the regular one-sided transposition shuffle by its selection of $k$ elements for $k$ consecutive permutations. 

Aside from the work done by Matheau-Raven on the one-sided transposition shuffle on the hyperoctahedral group \cite{M20}, Pang has computed the eigenvalues for card shuffling operators on the hyperoctahedral group \cite{MR4394678}. Pittet and Saloff-Coste studied random walks on the generalized symmetric groups and showed that they are closely related to random walks on simpler factor groups \cite{MR1905862}. Schoolfield, Jr. provided the rates of convergence for certain random walks on $G_{m,n}$ and independently derived the cutoff phenomenon for random walks on the hyperoctahedral group \cite{MR1922442}. 

Because we are unable to fully construct the operators which we have formulated in Section~\ref{sec:future-work} and show that they satisfy the conditions that would allow us to prove Conjectures~\ref{conjecture:lifting-eigenvectors_random} and \ref{conjecture:biasedOST}, it would be meaningful to complete the construction and prove that they do, in fact, suggest that the eigenvalues for the biased OST can be described by the technique of lifting eigenvalues. It would be interesting to develop a general framework for understanding how we can compute eigenvalues for different shuffling operators, other than the unbiased and biased OST shuffles, for the generalized symmetric group $G_{m,n}$.

\subsection{Organization of the Paper}
We now provide a brief overview of the contents of
the paper. The background material for our paper is contained in Section \ref{section:prelim}. It includes the necessary notation, definitions, and theorems to understand the unbiased OST shuffle and the generalized symmetric group. In Sections \ref{section:lower_bound} and \ref{section:upper_bound}, we prove the lower bound and upper bound for the cutoff for the OST shuffle on $G_{m,n}$ respectively. Section \ref{section:results} contains the formal statement of our main theorem, conjectured by Matheau-Raven. In Section \ref{section:results}, we prove this theorem by combining the bounds shown in Sections \ref{section:lower_bound} and \ref{section:upper_bound}. Sections~\ref{section:prelim}-\ref{section:results} cover our detailed proof for the cutoff time for the unbiased OST. In Section~\ref{sec:branchingrulesforgensymmgrp}, we lay the foundation for our attempt to tackle Conjecture~\ref{conjecture:lifting-eigenvectors_random} and $\ref{conjecture:biasedOST}$ and provide the preliminaries on the branching rules for the generalized symmetric group as well as the structure of its modules. Most importantly, we formally state Theorem~\ref{conjecture:informal_branchingrules} and prove it in Section~\ref{sec:branchingrulesforgensymmgrp}. Lastly, in Section~\ref{sec:future-work}, we formalize the concepts and operators as well as prove certain properties that they possess to set up the foundation for future work and attempts to prove Conjectures~\ref{conjecture:lifting-eigenvectors_random} and $\ref{conjecture:biasedOST}$.

\subsection{Acknowledgements}
The authors are grateful to Elizabeth Milićević for helpful discussions on the underlying representation theory, and for useful comments on an earlier version of this work. 

\section{Preliminaries}\label{section:prelim}
In this section, we introduce definitions and lemmas that are fundamental building blocks of the proof of Theorem \ref{thm:main}.


\subsection{Random Walks on Finite Group} In this subsection, we provide important notation, definitions and theorems for random walks on finite groups.
For the subsequent sections of this paper, we use $G$ to denote an arbitrary finite group.

\begin{definition}[{\cite[Total variation]{MR0841111}}]\label{def:variationdistance}
Denote $P^t$ the probability distribution of $X^t$ (resp. $(X^t)_{t \geq 0}$) the Markov chain on $G$, given that $X^0 = e$.
For any two distributions $P, Q$ on $G$, their \emph{total variation} distance is defined to be
\[\norm{P - Q}_{TV} := \frac{1}{2} \sum_{g \in G} \abs{P(g) - Q(g)}.\]
\end{definition}
Since we are interested in the time it takes to converge to the uniform distribution, the total variation distance that we are interested in is as follows:
\[d_{TV}(t):= \norm{P^t - U}_{TV}.\]

Now, we state a useful lemma for the total variation distance.

\begin{lemma}[{\cite[Lemma 7.10]{MR2466937}}]\label{lem:projection_tv_distance}
    Let $P$ and $Q$ be probability distributions on a finite group $G$. Let $f : G \to G'$
be a function on $G$ , where $G'$ is a relevant partition of $G$. Then,
\[\norm{P - Q}_{TV} \geq \norm{Pf^{-1} - Qf^{-1}}_{TV}.\]
\end{lemma}

In other words, the total variation distance between probability distributions can only decrease under projections.

\begin{definition}[{\cite[Separation Distance]{MR3530321}}]\label{def:separationdistance}
The worst-case \emph{separation distance} from stationary at time $t$ is defined as
\[d_{Sep}(t) := 1 - \min_{g, h \in G} P^t(hg^{-1})/U(h) = \norm{P^t - U}_{Sep}.\]
\end{definition}

\begin{theorem}[{\cite[Theorem 4.9]{MR2466937}}, the Convergence Theorem]\label{ref:convergence}
Suppose $P$ is irreducible and aperiodic with stationary distribution $U$.
Then there exists a constant $\alpha \in (0,1)$ and $C>0$ such that
\[\norm{P^t - U}_{TV} \leq C \alpha^t.\]
\end{theorem}

The probabilities in aperiodic and irreducible Markov chains always converge to some fixed probability distribution. By Matheau-Raven's work in \cite{M20}, the Markov chain of our interest demonstrates convergence.

Next, we introduce the \emph{cutoff phenomenon}. It is observed for some Markov chains, the total variation distance $\norm{P^t - U}_{TV}$ stays at its maximum 1, for a while, and then suddenly drops to $0$ as it converges.
The sudden drop from 1 to 0 is called the cutoff phenomenon.

\begin{definition}[{\cite[Mixing time]{MR2466937}}]\label{def:mixingtime}
The \emph{mixing time} with respect to a distance metric $d(t)$ is defined by
\[t_{mix}(\epsilon) := \min\{t: d(t) < \epsilon\}.\]
\end{definition}

\begin{definition}[{\cite[Cutoff]{MR2466937}}]\label{def:cutoff}
Suppose for a sequence of Markov chains indexed by natural numbers $n=1,2,\dots$, where $n$ represents the number of cards, the mixing time for the $n^{th}$ chain is denoted by $t^{(n)}_{mix}(\epsilon)$.
This sequence of chains has a \emph{cutoff} with window size $\alpha w_n$ for some $\alpha$ and $\lim_{n \to \infty}w_n/t_{mix}^{(n)} = 0$ when
\[\lim_{\alpha \to -\infty} \lim_{n \to \infty} d_n(t_{mix}^{(n)} + \alpha w_n ) = 1,\]
\[\lim_{\alpha \to \infty} \lim_{n \to \infty} d_n(t_{mix}^{(n)} + \alpha w_n ) = 0.\]
\end{definition}

\subsection{Generalized Symmetric Groups}\label{section:gen_sym_grp}
In this subsection, we provide the definition and examples of generalized symmetric groups. Readers should note we fix $m$ and $n$ to be natural numbers from Section \ref{section:gen_sym_grp} onwards as part of the definition for the generalized symmetric group. We also define $\xi$ to be the $m^{th}$ root of unity and denote the set of integers $\{1, \ldots, n\}$ as $[n]$.
\begin{definition}[{\cite[Section 1]{MR1422524}}, Generalized Symmetric Groups]\label{def:gen_symmetric_grp}

The generalized symmetric group $G_{m,n}$ is the group with order $m^nn!$ of all bijections $\varphi$ on $\left\{\xi^k i \mid i \in [n], k \in \Z_m\right\}$ such that $\varphi\left(\xi^k i\right) = \xi^k \varphi(i)$. In other words, $G_{m,n}$ is the set:
\[\left\{\left(\xi^{k_1}, \ldots, \xi^{k_n}, \sigma\right) \mid k_i \in \Z_m, \sigma \in S_n\right\}\]
with the operation:
\[\left(\xi^{k_1}, \ldots, \xi^{k_n}, \sigma\right)\left(\xi^{k_1'}, \ldots, \xi^{k_n'}, \sigma'\right) = \left(\xi^{k_1}\xi^{k_{\sigma(1)}'}, \ldots, \xi^{k_n}\xi^{k_{\sigma(n)}'}, \sigma\sigma'\right).\]
For any $\varphi = (\xi^{k_1}, \ldots, \xi^{k_n}, \sigma) \in G_{m,n}$ where $\sigma \in S_n$, we can define its image on $[n]$ in array notation as such:
\[\begin{bmatrix}1 & 2 & \ldots & n \\ \xi^{k_1}\sigma(1) & \xi^{k_2}\sigma(2) & \ldots & \xi^{k_n}\sigma(n)
\end{bmatrix}.\]

\end{definition}

The following are some examples of $G_{m, n}$ for $m = 1$ and $m = 2$:

\begin{example}$G_{1, n}$ is isomorphic to $S_n$.
\end{example}

\begin{example}$G_{2, n}$ is isomorphic to the hyperoctahedral group $B_n$.
\end{example}

Another interpretation of $B_n$ is the arrangement of a deck of $n$ cards where we now distinguish between cards that are facing up or down. 
The bijection $\varphi$ not only tells us where we should place the $i$\textsuperscript{th} card but also informs us whether the card should be facing up or down. For the case of $G_{m,n}$, the cards would have $m$ possible orientations and the bijection would tell us the number of rotations we should apply to the orientation of a particular card.

\subsection{One-sided Transposition Shuffle (OST)}
In this subsection, we introduce the OST shuffle, as well as its cutoff time on symmetric groups.
\begin{definition}[{\cite[Definition 1]{MR4312845}}, One-sided Transposition Shuffle] \label{def:ost}
The one-sided transposition (OST) shuffle $P_n$ is the random walk on $S_n$ generated by the following probability distribution on the conjugacy class of transpositions:
\[P_n(\sigma) = \begin{cases}\frac{1}{n}\cdot \frac{1}{j} & \text{ if $\sigma = (ij) \in S_n$ for some $ 1 \leq i \leq j \leq n$}, \\ 0 & \text{ otherwise.} \end{cases}\]
\end{definition}

Another way to interpret the OST shuffle is to first choose one of the $n$ cards in a deck. Suppose that card is the $j^{th}$ card from the top of the deck, we then choose our second card for the transposition to be one of the first $j$ cards (inclusive of the $j^{th}$ card). Since we are only looking at the cards above the $j^{th}$ card, it is an OST.

We use the convention that all transpositions with the form of $(ii)$ are equal to the identity $id$, and
therefore $P_n(id) = \frac{1}{n}\left(1 + \frac{1}{2} + \ldots + \frac{1}{n}\right) = H_n/n$, where $H_n$ denotes the $n^{th}$ harmonic number. 

\begin{theorem}[{\cite[Theorem 3]{MR4312845}}, Cutoff for OST on $S_n$] \label{thm:cutoff_ost_symmetric_group}
The OST shuffle $P_n$ on $S_n$ exhibits a cutoff at time $n \log(n)$. For any $c_1 > 0, c_2 > 2$,
\[\limsup_{n\to \infty}{\norm{P_n^{n\log(n) + c_1 n}- U_n}}_{TV} \leq \sqrt{2}e^{-c_1},\]
\[\text{and }\,\liminf_{n\to\infty}{\norm{P_n^{n\log(n) - n\log\log n - c_2 n}- U_n}}_{TV} \geq 1 - \frac{\pi^2}{6(c_2 - 2)^2}.\]
\end{theorem}

\subsection{One-sided Transposition Shuffle on the Generalized Symmetric Group}
In this subsection, we provide the formal definition of the OST on the generalized symmetric group $G_{m,n}$ based on Definitions \ref{def:gen_symmetric_grp} and \ref{def:ost}.

\begin{definition}[OST on $G_{m,n}$]\label{def:ost_gensymgroup}
For any $(ij) \in S_n$ where $1 \leq i \leq j \leq n$, define $G_{m,n}^{(ij)}$ to be the following subset of $G_{m,n}$:
\[G_{m,n}^{(ij)} = \left\{\left(\underbrace{1, \ldots, 1}_{i-1 \text{ terms}}, \xi^{k}, \underbrace{1, \ldots, 1}_{j-i-1 \text{ terms}}, \xi^{k}, \underbrace{1, \ldots, 1}_{n-j \text{ terms}}, (ij)\right) \,\middle\vert\, k \in \mathbb{Z}_m\right\}.\]
Let $\sigma$ be an arbitrary element in $G_{m,n}$. The OST shuffle $OST_{m,n}$ is the random walk on $G_{m,n}$ generated by the following probability distribution:
\[OST_{m,n}(\sigma) = \begin{cases}\frac{1}{n\cdot j\cdot m} & \,\sigma \in G_{m,n}^{(ij)} \text{ for some } (ij) \in S_n\\
0 &\, \text{ otherwise.}\end{cases}\]

Another way to interpret the above definition is to first perform the standard OST shuffle on the deck of $n$ cards as per Definition \ref{def:ost}. Subsequently, uniformly pick $k$ from $\Z_m$ and flip the transposed cards $k$ times. Note that the cards in the deck have $m$ possible orientations. Suppose the indexing for the different orientations starts at $0$. Then, if the card started off at orientation number 2, the flipping will send it to orientation number $(2 + k) \bmod{m}$.
\end{definition}

\section{Lower Bound for OST Shuffle on The Generalized Symmetric Group.}\label{section:lower_bound}
In this section, we begin by stating and proving a lemma on the relationship between $OST_{m,n}$ and $OST_{1,n}$ before finding a lower bound for $OST_{m,n}$ using the lower bound for $OST_{1,n}$.

\begin{lemma}\label{lem:projection_lemma}
Given an arbitrary $\eta \in S_n$, define an onto homomorphism $\psi: G_{m,n} \to S_n$ such that $\psi\left(\left(\xi^{k_1}, \ldots, \xi^{k_n}, \eta\right)\right) = \eta$ for $k_1, \ldots, k_n \in \mathbb{Z}_m$. Let $OST_{m, n}$ be the OST shuffle on $G_{m, n}$. Define $t$ to be the number of convolutions of the $OST_{m,n}$ probability distribution.
Then, for all $t \geq 1$, we have
\[\sum_{\sigma \in \psi^{-1}(\eta)}OST^t_{m, n}(\sigma) = OST^t_{1, n}(\eta).\]

\begin{proof}
Let us prove the above lemma by inducting on $t$. For the base case, we have $t = 1$. For any $\eta = (ij) \in S_n$ where $1 \leq i \leq j \leq n$, we have the following:
\[\psi^{-1}(\eta) = \bigcup_{k_i, k_j \in \Z_m}\left\{\left(\underbrace{1, \ldots, 1}_{i-1 \text{ terms}}, \xi^{k_i}, \underbrace{1, \ldots, 1}_{j-i-1 \text{ terms}}, \xi^{k_j}, \underbrace{1, \ldots, 1}_{n-j \text{ terms}}, \eta\right)\right\}.\]
We know from Definition \ref{def:ost_gensymgroup} that the only $m$ elements in $G_{m,n}$ with non-zero probability belong to $G_{m,n}^{(ij)}$. In fact, we know that all $m$ of them have the same probability of $\frac{1}{n\cdot m \cdot j}$. Then, 
\begin{align*}
    \sum_{\sigma \in \psi^{-1}(\eta)}OST_{m, n}(\sigma) &= m \cdot \left(\frac{1}{n\cdot m \cdot j}\right) \\
    &= OST_{1, n}(\eta).
\end{align*}
Having shown that the base case holds, let us now state the induction hypothesis: suppose $\sum_{\sigma \in \psi^{-1}(\eta)}OST^k_{m, n}(\sigma) = OST^k_{1, n}(\eta)$ for some $k \geq 1$. For the inductive step, we know from the convolution of the $OST_{1,n}$ probability distribution that 

\begin{align*}
    OST_{1,n}^{k+1}(\eta) &= \sum_{g\in G_{1,n}}OST_{1,n}(\eta g^{-1}) OST_{1,n}^k(g) \\
    &= \sum_{g \in G_{1,n}}\left(\sum_{\tau \in \psi^{-1}(\eta g^{-1})}OST_{m,n}(\tau)\right)\left(\sum_{\tau' \in \psi^{-1}(g)}OST^k_{m,n}(\tau')\right) \\
    &= \sum_{g \in G_{1,n}}\sum_{\tau \in \psi^{-1}(\eta g^{-1})}\sum_{\tau' \in \psi^{-1}(g)}OST_{m,n}(\tau)OST^k_{m,n}(\tau') \\
    &= \sum_{\sigma \in \psi^{-1}(\eta)}OST_{m,n}^{k+1}(\sigma)
\end{align*}
where the second equality comes from the base case and the induction hypothesis, and the last equality comes from using convolution and considering the possible walks to the preimage of the onto homomorphism $\psi$ for $\eta$.
\end{proof}
\end{lemma}
By adapting Lemma 4.3.3 in \cite{M20}, we can now prove a proposition for the lower bound for the OST shuffle on $G_{m,n}$. For the sake of notation brevity, we also abbreviate the uniform distribution on $G_{m,n}$ as $U_{m,n}$.
\begin{proposition}\label{lem:lowerbound_OST_gensymgroup}
The OST shuffle $OST_{m,n}$ on $G_{m,n}$ satisfies the following lower bound for any scaling factor $c > 2$:
\[\liminf_{n\to\infty}\norm{OST_{m,n}^{n\log n - n\log\log n - c n} - U_{m,n} }_{TV} \geq 1 - \frac{\pi^2}{6(c-2)^2}.\]
\begin{proof}
By Lemma \ref{lem:projection_tv_distance}, we know that the total variation distance can only decrease under projections. Hence, let us define a surjective homomorphism $\psi: G_{m, n} \to S_n$ by ignoring the first $n$ elements in the tuple for $G_{m, n}$, i.e., for $\sigma \in G_{m, n}$, $\psi((\xi^{k_1}, \ldots, \xi^{k_n}, \eta)) = \eta$ for any $k_1, \ldots, k_n \in \Z_m$ and $\eta \in S_n$. Given an arbitrary $\eta \in S_n$, consider the preimage of $\psi$ and $t \geq 1$, we have the following from Lemma \ref{lem:projection_lemma}:
\[\sum_{\sigma \in \psi^{-1}(\eta)}OST^t_{m,n}(\sigma) = OST^t_{1, n}(\eta).\]
Hence, we can now conclude the following:
\begin{align*}
    \norm{OST_{m,n}^t - U_{m,n}}_{TV} &= \frac{1}{2}\sum_{\sigma \in G_{m,n}}\left|OST_{m,n}^t(\sigma) - U_{m,n}(\sigma)\right|\\
    &\geq \frac{1}{2}\sum_{\eta \in S_n}\left|\sum_{\sigma \in \psi^{-1}(\eta)}OST_{m,n}^t(\sigma) - U_{m,n}(\sigma)\right| \\
    &= \frac{1}{2}\sum_{\eta \in S_n}|OST_{1,n}(\eta) - U_{1,n}| \\
    &= \norm{OST_{1,n}^t - U_{1,n}}_{TV} 
\end{align*}
where the first equality follows from Definition \ref{def:variationdistance}, the second equality follows from Lemma \ref{lem:projection_lemma}, and the last equality follows from Definition \ref{def:variationdistance}. With the above inequality, the lower bound in the proposition statement now follows directly from Theorem \ref{thm:cutoff_ost_symmetric_group}.
\end{proof}
\end{proposition}
\section{Upper Bound for OST Shuffle on The Generalized Symmetric Group.}\label{section:upper_bound}
In this section, we find the upper bound by providing a connection between OST on generalized symmetric group and the coupon collector's problem.
\begin{definition}
Let $(X^t)_{t \in \N}$ be a Markov chain on the generalized symmetric group $G_{m,n}$ driven by the OST shuffle $OST_{m,n}$. Define $(Y^t)_{t \in \N}$ to be a Markov chain on $G_{m,n}$ such that $Y^t$ is the inverse of $(X^t)$ for all $t$. 
\end{definition}

\begin{definition}[{\cite[Stopping time]{MR0964069}}]
A \emph{stopping time} is a function $T:(X^t)_{t\in\N} \to \N$ such that if $T((X^t)) = j$ then $T((\tilde{X}^t)) = j$ for all $(\tilde{X}^t)$ with $\tilde{X}_i = X_i$ for $1 \leq i \leq j$.
\end{definition}
Intuitively, a stopping time looks at a Markov chain and decides if a condition has been met based on its first $j$ steps, without looking at future steps.
\begin{definition}[{\cite[Strong stationary time]{MR0841111}}]\label{def:MR0841111}
A \emph{strong stationary time} for a Markov chain is a stopping time $T$ such that $X^T$ is stationary and independent of $T$, i.e., for any $g \in G,$
\[\mathbb{P}(X^k = g \mid T = k) = U(g).\]
\end{definition}

We also state a useful lemma that uses strong stationary times to bound the separation distance and total variation distance from above.

\begin{lemma}[{\cite[Lemma 1.1.31]{M20}}]\label{lem:stationary_time_bound_sep_dist}
Let $T$ be a strong stationary time for an irreducible, aperiodic random walk on $G$ with a probability distribution $P$. Let $U$ be the uniform distribution on $G$. The following holds for all $t \geq 1$:
\[\norm{P^t - U}_{TV} \leq \norm{P^t - U}_{Sep} \leq \mathbb{P}\left(T > t\right)\]
\end{lemma}

Recall from Definition \ref{def:gen_symmetric_grp} that $G_{m,n}$ is the set of bijections $\varphi$ on $\{\xi^k i \mid k \in \Z_m, i \in [n]\}$ such that $\varphi(\xi^k i) = \xi^k \varphi(i)$ where $\xi$ is the $m^{th}$ root of unity. 
Similar to how a strong stationary time was constructed in \cite{M20} for the random walk on $S_n$, in order to construct a strong stationary time for our random walk on $G_{m,n}$ with the OST shuffle, we need to condition on exact knowledge of the positions above position $j$ in the card deck at time $t$, i.e., the random variables $Y^t(\xi^k i)$ for $k \in \Z_m$ and $j < i \leq n$ where $Y^t(\xi^k i) = \xi^k Y^t(i)$.

Adapting Definition 3.7.3 from \cite{M20}, we propose a property that tells us whether the known information about the deck strictly above position $j$ at time $t$ is equally likely to be any of the remaining cards.
\begin{definition}[Property $\mathcal{P}_j$]\label{def:property_j}
Consider the random variables $Y^t(\xi^k i)$ for $k \in \Z_m$ and $j < i \leq n$ where $Y^t(\xi^k i) = \xi^k Y^t(i)$. Define the set of cards which we know can be in positions 1 to $j$ at time $t$ as\[\mathcal{G}_j^t = \{\xi^k i \mid i \in [n], k \in \Z_m\} \setminus \{Y^t(\xi^k i) \mid j < i \leq n, k \in \Z_m\}.\]  
Then, $(Y^t)_{t \in \N}$ satisfies property $\mathcal{P}_j$ at time $t$ if:
\[\mathbb{P}(Y^t(j) = x \mid Y^t(i) \text{ for all }j < i \leq n) = \begin{cases}\frac{1}{mj} &\,x \in \mathcal{G}^t_{j}, \\ 0 &\,\text{otherwise.}\end{cases}\]
\end{definition}
Having defined property $\mathcal{P}_j$ for $G_{m,n}$, we now prove a strengthened version of Lemma 4.4.12 from \cite{M20} for $G_{m,n}$.
\begin{proposition} \label{prop:satisfy_property_pj}
Let $T_j$ be the first time we select from position $j$ on the first draw when performing the OST shuffle on $G_{m,n}$. If $T_j \leq t$, then $(Y^t)_{t \in \N}$ satisfies property $\mathcal{P}_j$ at time $t$.
\begin{proof}
We proceed by induction on $t$, i.e., the Markov chain $(Y^t)_{t \in \N}$ satisfies property $\mathcal{P}_j$ for all times after $t$ once it holds for some $t \geq 0$.

Consider the time $T_j$ where we must have applied a transposition $\xi^{k}_i\xi^{k'}_j(ij)$ as part of our Markov chain with $i \leq j$ and $k, k' \in \Z_m$. 
The probability of picking any one of the transpositions $\tau^{T_j} = \xi^{k}_i\xi^{k'}_j(ij)$ at time $T_j$ for all $i \leq j$ is
\[\mathbb{P}(\tau^{T_j} = \xi^{k}_i\xi^{k'}_j(ij)) = \frac{1}{mj}.\]

Recall from Definition \ref{def:property_j} that $\mathcal{G}^t_j$ denotes the set of cards which we know can be in positions 1 to $j$ at time $t$. 
Hence, the card in position $j$ at time $T_j$ has equal likelihood of being any of the cards in $\mathcal{G}^{T_{j-1}}_j = \mathcal{G}^{T_{j}}_j$ (which contains cards including their $m$ possible orientations). 
Thus, it becomes clear that:
\[\mathbb{P}(Y^{T_j}(j) = x \mid Y^{T_j}(i) \text{ for all } j < i \leq n) = \begin{cases} \frac{1}{mj} &\,\text{ if }x \in \mathcal{G}^{T_j}_{j} \\ 0 &\,\text{ otherwise,}\end{cases}\]
which implies that $\mathcal{P}_j$ holds at time $T_j$.

Our inductive hypothesis allows us to assume that the Markov chain satisfies property $\mathcal{P}_j$ at some arbitrary time $t \geq 0$. 
Then,
\begin{align*}
    &\mathbb{P}(Y^{t+1}(j) = x \mid Y^{t+1}(i) \text{ for all } j < i \leq n) \\
    = \sum_{\sigma \in G_{m,n}} &\mathbb{P}(\tau^{t+1} = \sigma) \mathbb{P}(Y^{t+1}(j) = x \mid Y^{t+1}(i) \text{ for all } j < i \leq n, \tau^{t+1} = \sigma).
\end{align*}

We now show in inductive steps that the property $\mathcal{P}_j$ still holds at time $t+1$ by splitting our analysis into cases depending on the permutation $\tau^{t+1}$.
We split the cases by the positions of $a$ and $b$ relative to $j$ as follows:
\begin{align*}
    \text{case 1: }& \tau^{t +1} \in \{(a j), \xi_a^k \xi_j^{k'} (aj) \mid a \leq j\};\\
    \text{case 2: }& \tau^{t+1} \in \{(ab), \xi_a^k \xi_b^{k'} (ab) \mid a,b < j\};\\
    \text{case 3: }& \tau^{t+1} \in \{(ab), \xi_a^k \xi_b^{k'} (ab) \mid j < a,b\};\\
    \text{case 4: }& \tau^{t+1} \in \{(ab), \xi_a^k \xi_b^{k'} (ab) \mid a \leq j < b\}.
\end{align*}

For case 1, this random walk satisfies property $\mathcal{P}_j$ for the same reason as the $T_j$ case.
For case 2, we know that the position of the cards above $j$ did not change.
Hence, $Y^t(i) = Y^{t+1}(i)$ for all $j < i \leq n$.
For case 3, 
\begin{align*}
    Y^t(\xi^k a) = Y^{t+1} (\xi^k b), & \text{ } Y^t(\xi^{k'} a) = Y^{t+1} (\xi^{k'} b),\\
    Y^t(\xi^k b) = Y^{t+1} (\xi^k a), & \text{ } Y^t(\xi^{k'} b) = Y^{t+1} (\xi^{k'} a),
\end{align*}
with $Y^{t}(i) = Y^{t+1}(i)$ for all $j < i \leq n$.

Hence, we have
\[\mathbb{P}(Y^{t+1}(j) = x \mid Y^{t+1}(i) \text{ for all } j < i \leq n, \tau^{t+1} \in \{(ab) \xi_a^k \xi_b^{k'} (ab) \mid a, b < j \text{ or both } a, b > j\}),\]
which is equivalent to
\[\mathbb{P}(Y^t(j) = x \mid Y^t(i) \text{ for all } j < i \leq n) = \begin{cases} \frac{1}{mj} &\,\text{ if }x \in \mathcal{G}^{T_j}_{j}, \\ 0 &\,\text{ otherwise.}\end{cases}\]

For case 4, we cannot yield the relationship between $Y^t(i)$ and $Y^{t+1}(i)$ without further assumptions.

Keep $b > j$ fixed and consider a card $C$ in $\mathcal{G}_j^{t+1}$.
Suppose that $\tau^{t+1}$ is the action that moves $C$ from position $b$ to another position below $j$.
Since we know that $X^t = \tau^{t+1} X^{t+1}$, we have the following:
\[\tau^{t+1} = \mqty(X^{t+1}(C) &b).\]

Let $C$ iterate through all cards in $\mathcal{G}_j^{t+1}$.

Then, we know that $Y^t(b) = C$.
Hence, all the cards above $j$ except $C$ did not move position. 
Then,
\[Y^t(i) = Y^{t+1}(i) \text{ for } j 
< i \leq n, i \neq b.\]
Therefore,
\[\mathcal{G}_j^t = (\mathcal{G}_j^{t+1} \sqcup \{Y^{t+1}(\xi_b^k)\})\setminus\{\xi_b^k C\}.\]

Now, consider the probability
\[\mathbb{P}\left(Y^{t+1}(j) = x \mid Y^t(i) \text{ for all } j < i \leq n, \tau^{t+1} = \mqty(X^{t+1}(C) &b)\right).\]

When $x = \xi_b^k C$, then the card $C$ has to be in $\xi_b^k b$ at time $t$.
Hence, $Y^{t+1}(b) \in \mathcal{G}_j^t$; when $x \in \mathcal{G}_j^{t+1} \setminus 
\{\xi_b^k C\}$, then we know that $x$ did not change position from time $t$ to $t+1$, so $x \in \mathcal{G}_j^t$.

Therefore,
\[\mathbb{P}\left(Y^{t+1}(j) = x \mid Y^t(i) \text{ for all } j < i \leq n, \tau^{t+1} = \mqty(X^{t+1}(C) &b)\right) = \begin{cases}1/mj & x \in \mathcal{G}_j^t,\\ 0 & \text{otherwise}.\end{cases}\]

Since we can let $C$ iterate through all cards in $\mathcal{G}_j^{t+1}$, we know that for all $a < j \leq b$, we have shown that:
\[\mathbb{P}\left(Y^{t+1}(j) = x \mid Y^t(i) \text{ for all } j < i \leq n, \tau^{t+1} = (ab)\right) = \begin{cases}1/mj & x \in \mathcal{G}_j^t,\\
0 & \text{otherwise}.\end{cases}\]
Hence, we have proven the last case, and combined with all three cases above, we find that:
\[\mathbb{P}\left(Y^{t+1}(j) = x \mid Y^t(i) \text{ for all } j < i \leq n \right) = \begin{cases}1/mj & x \in \mathcal{G}_j^t,\\
0 & \text{otherwise}.\end{cases}\]

This concludes the inductive step.
\end{proof}
\end{proposition}

With the above proposition, we can now state that \[T = \min\left\{t \geq 0\,:\, t \geq T_j \text{ for all }j\right\}\]
is a stopping time for $(Y^t)_{t \in \N}$. Before we proceed to bound the separation distance mixing time of the OST shuffle on $G_{m,n}$ from above, we first state a useful proposition from the coupon collector's problem.

\begin{proposition}[{\cite[Propositions 2.3 and 2.4]{MR2466937}}]\label{prop:coupon_collector}
    Consider a collector attempting to collect a complete set of coupons. Assume that each new coupon is chosen uniformly and independently from the set of n possible types, and let $\ell$ be the (random) number of coupons collected when the set first contains every type. For any $c > 0$.
    \[\mathbb{P}\left(\ell >  n\log n + cn\right) \leq e^{-c}\]
\end{proposition}

Now, we are ready to prove the following proposition that states the upper bound for the separation distance mixing time of the OST shuffle on $G_{m,n}$. Some parts of our proof were inspired by the proof sketch in Section 4.5 of \cite{M20}.

\begin{proposition}\label{lem:upperbound_OST_gensymgroup}
The OST shuffle $OST_{m,n}$ on $G_{m,n}$ satisfies the following upper bound for any $c > 0$:
\[\limsup_{n\to\infty}\norm{OST_{m,n}^{n\log n + cn} - U_{m,n} }_{Sep} \leq e^{-c}.\]
\begin{proof}
Since we now know that $T = \min\{t \geq 0\,:\, t \geq T_j \text{ for all }j\}$
is a stopping time for $(Y^t)_{t \in \N}$, we know that $(Y^T)$ satisfies property $\mathcal{P}_j$ for all $j$ by Proposition \ref{prop:satisfy_property_pj}. 

For any $\sigma \in G_{m,n}$,
\begin{align*}
    \mathbb{P}(Y^t = \sigma^{-1} \mid t \geq T) &= \mathbb{P}(\cap_{j=1}^n\{Y^t(j) = \sigma^{-1}(j)\} \mid t \geq T) \\
    &= \prod_{j=1}^n\mathbb{P}(Y^t(j) = \sigma^{-1}(j) \mid \cap_{i = j+1}^n\{Y^t(i) = \sigma^{-1}(i)\}, t \geq T) \\
    &= \prod_{j=1}^n \frac{1}{mj} = \frac{1}{m^nn!} = U_{m,n}(\sigma).
\end{align*}
This implies that $T$ is a strong stationary time for $(Y^t)_{t \in \N}$. By construction of $(Y^t)_{t \in \N}$,  $T$ is also a strong stationary time for $(X^t)_{t \in \N}$.

Recalling that an OST shuffle on $G_{m,n}$ involves drawing two cards in sequential order, $T$ is defined to be the first time that we have chosen every possible position $j$ for our first draw during the OST shuffle on $G_{m,n}$. By Definition \ref{def:ost_gensymgroup}, our first draw chooses position $j$ uniformly on $[n]$. Thus, $T$ is similar to the $\ell$ stated in Proposition \ref{prop:coupon_collector} which gives us $\mathbb{P}(T > n\log n + cn) \leq e^{-c}$. Substituting $OST_{m,n}^{n\log n + cn}$ and $T$ into Lemma \ref{lem:stationary_time_bound_sep_dist} gives us 
\[\norm{OST_{m,n}^{n\log n + cn} - U_{m,n}}_{Sep} \leq \mathbb{P}(T > n\log n + cn).\]
Now, we can conclude that 
\[\limsup_{n\to \infty}\norm{OST_{m,n}^{n\log n + cn} - U_{m,n}}_{Sep} \leq \limsup_{n\to\infty}\mathbb{P}(T > n\log n + cn) \leq e^{-c}.\]
\end{proof}
\end{proposition}
Having proven the upper bound for the separation distance mixing time of the OST shuffle on $G_{m,n}$, we are now ready to prove our main theorem in Section \ref{section:results}.

\section{Cutoff for OST Shuffle on the Generalized Symmetric Group.}\label{section:results}
In this section, we state and prove the formal version of Theorem \ref{thm:main}. Our approach utilizes Propositions \ref{lem:lowerbound_OST_gensymgroup} and \ref{lem:upperbound_OST_gensymgroup} proved in Sections \ref{section:lower_bound} and \ref{section:upper_bound} for the lower and upper bounds for the cutoff time.

\begin{theorem}[{\cite[Conjectured in Section 4.5]{M20}}]\label{thm:main_formal}
The unbiased OST shuffle on the generalized symmetric group $G_{m,n}$ exhibits a cutoff in total variation distance and separation distance at time $n\log(n)$.
\end{theorem}

\begin{proof}
For arbitrary constants $c_1>0$ and $c_2>2$, we define two functions
\[F_1(n,c_1) = n\log n + c_1n,\]
\[F_2(n, c_2) = n\log n - n\log \log n - c_2n.\]
By Lemma \ref{lem:stationary_time_bound_sep_dist}, we know that the separation distance is an upper bound for the total variation distance. From Propositions \ref{lem:lowerbound_OST_gensymgroup} and \ref{lem:upperbound_OST_gensymgroup}, we have the following :
\begin{align*}\lim_{n\to\infty}\norm{OST_{m,n}^{F_1(n, c_1)} - U_{m,n}}_{Sep}&\leq\lim_{n\to\infty}\norm{OST_{m,n}^{F_1(n, c_1)} - U_{m,n}}_{Sep}\leq e^{-c_1}, \\
\lim_{n\to\infty}\norm{OST_{m,n}^{F_2(n, c_2)} - U_{m,n}}_{Sep} &\geq\lim_{n\to\infty}\norm{OST_{m,n}^{F_2(n, c_2)} - U_{m,n}}_{TV}\geq 1 - \frac{\pi^2}{6(c_2-2)^2}.
\end{align*}
Since $\lim_{c_1\to\infty}e^{-c_1} = 0$ and $\lim_{c_2\to\infty}1 - \frac{\pi^2}{6(c_2-2)^2} = 1$, by Definition \ref{def:cutoff}, we have proven that the OST shuffle on $G_{m,n}$ exhibits a cutoff in total variation distance and separation distance at time $n \log n$. 
\end{proof}

\section{Branching Rules for the Generalized Symmetric Group}\label{sec:branchingrulesforgensymmgrp}
In this section, we prove the branching rules for the specht modules of the generalized symmetric group (Theorem~\ref{thm:branching_rules_gensymgroup}) by first laying down the preliminaries for the structure of modules for the generalized symmetric group before using the branching rules for the generalized symmetric group to prove the theorem.
\subsection{Structure of Modules for the Generalized Symmetric Group}
In this subsection, we provide definitions and notations needed to study the structure of modules for the generalized symmetric group $G_{m,n}$. We first introduce partitions before delving into the specht modules for the generalized symmetric group and their branching rules. As a convention,
throughout the rest of this paper, we assume that $\xi$ is a primitive m-th root of unity
\subsubsection{Partitions for Generalized Symmetric Group}
In this subsection, we first state the definition of $m$-partitions and define the dominance ordering on $m$-partitions.
\begin{definition}[{\cite[Definition 1.2]{MR1422524}}, $m$-partition]\label{def:m_partition}
Let $m,n \in \N$. An $m$-partition of $n$, denoted by $\lambda^{[m]}$, is a tuple of partitions $\lambda^{[m]} = (\lambda^{(1)}, \ldots, \lambda^{(m)})$ such that $\sum_{k=1}^m |\lambda^{(k)}| = n$. In $\lambda^{[m]}$, let $\lambda^{(i)}$ be the $i^{th}$ constitutent of $\lambda^{[m]}$.
\end{definition}
In other words, the sum of the size of each partition in the tuple $\overline{\lambda}$ has to add up to $n$.
\begin{definition}[{\cite[Definition 2.1]{MR1422524}}, Dominance ordering on $m$-partitions]\label{def:dominance_m_partition} Let $\lambda^{[m]}, \mu^{[m]}$ be $m$-partitions. The dominance ordering on $m$-partitions is defined as:
\[\lambda^{[m]}\trianglerighteq \mu^{[m]} \iff \begin{cases} & |\lambda^{(1)}| > |\mu^{(1)}|, \\
\text{ or } & |\lambda^{(k)}| = |\mu^{(k)}| \text{ for } 1 \leq k \leq i, \text{ and } |\lambda^{(i+1)}| > |\mu^{(i+1)}| \text{ for some } i \in [m-2], \\ 
\text{ or } & |\lambda^{(k)}| = |\mu^{(k)}| \text{ and } \lambda^{(k)} \trianglerighteq \mu^{(k)} \text{ for all } k \in [m].
\end{cases}\]
\end{definition}
\subsubsection{Specht Modules for Generalized Symmetric Groups}
In this subsection, we first state the definition of polytabloids and then state two theorems that would be useful for proving Young's rule for the permutation and simple modules for $G_{m,n}$.
\begin{definition}[{\cite[Definition 3.1]{MR1422524}}, Polytabloid]
Let $t$ be a $\lambda^{[m]}$-tableau and $\varphi \in G_{m,n}$. Define $\kappa_t \in \mathbb{C} [G_{m,n}]$ by
\[\kappa_t = \sum_{\varphi \in C_t}\xi^{-f(\varphi)}(\operatorname{sgn}\varphi)\varphi\]

Given any $\varphi$ can be expressed as the product of disjoint cycles $\theta_1\ldots\theta_t$ where each of the $\theta_i$ can be expressed as such:
\[\begin{bmatrix}b_{i1} & b_{i2} & \ldots & b_{in} \\ \xi^{s_{i1}}b_{i2} & \xi^{s_{i2}}b_{i3} & \ldots & \xi^{s_{in}}b_{i1}
\end{bmatrix}\]
where $b_{ij} \in [n], s_{ij} \in [m]$, and $k_i$ is the length of the cycle.

The function $f$ is defined as such :
\begin{align}
    f(\varphi) &= \sum_{i=1}^t f(\theta_i) \bmod{m}\\
    &= \sum_{i=1}^t \sum_{j = 1}^{k_i}s_{ij} \bmod{m}
\end{align}
 
The $\lambda^{[m]}$-polytabloid $e_t$ associated with the tableau $t$ is given by $e_t = \kappa_{t}[t]$ where $[t]$ is the $\lambda^{[m]}$-tabloid that contains $t$ as per \cite{MR1422524}.
\end{definition}

Now, we state two theorems from \cite{MR1422524} that provide us with the homomorphism between the Sphect modules to the permutation modules and allow us to consider the set of Sphect modules and their relationship to $\mathbb{C}[G_{m,n}]-$modules. 
\begin{theorem}[{\cite[Proposition 3.20]{MR1422524}}]\label{thm:gensymgroup_schur}
Let $\lambda^{[m]}$, $\mu^{[m]}$ be $m$-partitions such that $|\lambda^{(i)}| = |\mu^{(i)}|$ for all $i \in [m-1]$. Suppose the field of scalars is $\mathbb{C}$ and $\psi:S^{\lambda^{[m]}}\to M^{\mu^{[m]}}$ is a non-zero homomorphism. Thus $\lambda^{[m]} \trianglerighteq \mu^{[m]}$ and if $\lambda^{[m]} = \mu^{[m]}$, $\psi$ is multiplication by a scalar.
\end{theorem}

\begin{theorem}[{\cite[Theorem 3.21]{MR1422524}}]\label{thm:basis_gensymgroup_specht}
The set of Specht modules $S^{\lambda^{[m]}}$ (for $m$-partitions $\lambda^{[m]}$) form a complete set of irreducible $\mathbb{C} [G_{m,n}]$-modules. 
\end{theorem}

Having defined $m$-partitions as well as the dominance ordering on $m$-partitions, let us now prove the $G_{m,n}$ version of Young's Rule.

\begin{lemma}\label{lem:young_rule_gen_sym_group}
The permutation and simple modules for $G_{m,n}$ are indexed by $m$-partitions of $n$. The permutation module and simple module corresponding to $\mu^{[m]}$ and $\lambda^{[m]}$ are denoted $M^{\lambda^{[m]}}$ and $S^{\lambda^{[m]}}$ respectively. Furthermore, the permutation and simple modules respect Young's rule, that is 
\[M^{\mu^{[m]}} \cong \bigoplus_{\lambda^{[m]} \trianglerighteq \mu^{[m]}} K_{\lambda^{[m]}, \mu^{[m]}} S^{\lambda^{[m]}},\]
for constants $K_{\lambda^{[m]}, \mu^{[m]}} \in \N \cup \{0\}$. If $\lambda^{[m]} = \mu^{[m]}$, $K_{\lambda^{[m]}, \mu^{[m]}} = 1$. 
\end{lemma}
\begin{proof}
By Theorem \ref{thm:basis_gensymgroup_specht}, we can decompose $M^{\mu^{[m]}}$ into the complete set of simple modules $S^{\lambda^{[m]}}$.

If $S^{\lambda^{[m]}}$ appears with a non-zero coefficient in the decomposition of $M^{\mu^{[m]}}$, then we naturally have a non-zero homomorphism $\psi:S^{\lambda^{[m]}}\to M^{\mu^{[m]}}$. By Theorem \ref{thm:gensymgroup_schur}, we know that if $S^{\lambda^{[m]}}$ appears as a summand of $M^{\mu^{[m]}}$, we are guaranteed that $\lambda^{[m]} \trianglerighteq \mu^{[m]}$. 
If $\lambda^{[m]} = \mu^{[m]}$, by Theorem \ref{thm:gensymgroup_schur}, we know that the homomorphism $\psi:S^{\lambda^{[m]}}\to M^{\mu^{[m]}}$ is multiplication by a scalar and hence there must be only be one copy of $S^{\lambda^{[m]}}$ in $M^{\mu^{[m]}}$. Hence, $K_{\lambda^{[m]}, \mu^{[m]}} = 1$ for $\lambda^{[m]} = \mu^{[m]}$.
\end{proof}

\subsubsection{Branching Rules for the Specht modules of the Generalized Symmetric Group}\label{sec:branchingrules}

In this subsection, we discuss the branching rules for the Specht modules of the generalized symmetric group. We do that by first defining the outer product of representations before stating the Littlewood-Richardson Rule for $G_{m,n}$. Finally, we prove the branching rules for the Specht modules of the generalized symmetric group.

Similar to the outer product introduced in \cite{MR0125885} for symmetric groups and \cite{MR0491917} for hyperoctahedral groups, let us define the outer product of $G_{m,n}$. Note that the outer product is associative and commutative.

\begin{definition}\label{def:outer_product}
Let $\lambda$ be a representation of $G_{m,k}$ and $\tau$ be a representation of $G_{m,n-k}$. Then, the outer product of $\lambda$ and $\tau$, $\operatorname{Ind}_{G_{m,k} \times G_{m,n-k}}^{G_{m,n}}\lambda \otimes \tau$, is a representation of $G_{m,n}$ where the outer product is denoted by $\lambda \# \tau$.
\end{definition}

Having stated the definition of outer products for group representations, we now state the Littlewood-Richardson Rule for the generalized symmetric group that would be integral in our effort to prove Theorem~\ref{thm:branching_rules_gensymgroup}.

\begin{theorem}[{\cite[Theorem 5]{MR0237669}}, Littlewood-Richardson Rule for $G_{m,n}$] \label{thm:littwoodrichardson}
Each diagram constructed according to the statements described below defines an irreducible component of $\lambda \# \tau$ and all components are obtained in this manner.
\begin{enumerate}
    \item To each tableau $\lambda^{(i)}$ of $\lambda^{[m]}$, add the symbols of the first row of a tableau $\tau^{(i)}$ of $\tau^{[m]}$ for $i = 1, 2, \ldots, m$. These may be added to one row or divided into any number of sets, preserving their order, the first set being added to one row of $\lambda^{(i)}$, the second set to a subsequent row, the third to a row subsequent to this, and so on. After the addition, no row of the compound tableau may contain more symbols than a preceding row, and no two added symbols may appear in the same column. Next, add the second row of $\tau^{(i)}$, according to the same rules, followed by the remaining rows in succession until all the symbols of $\tau^{[i]}$ have been used.
    \item These additions must be such that each symbol from $\tau^{(i)}$ shall appear in a later row of the compound tableau than that occupied by the symbol immediately above it in $\tau^{[i]}$, for $i = 1, 2, \ldots, m$.
\end{enumerate}
\end{theorem}

Now, we are ready to formally state the branching rules for the Specht modules for the generalized symmetric group and prove it.

\begin{theorem}[Branching Rules for the Sphect Modules of $G_{m,n}$]\label{thm:branching_rules_gensymgroup}
Let $n \geq 2$, $\lambda^{[m]} \vdash n$ be $m$-partitions of $n$, $\mu^{[m]} \vdash n + 1$ be $m$-partitions of $n+1$, and $\tau^{[m]} \vdash 1$ be $m$-partitions of 1. The branching rules for the Specht modules of $G_{m,n}$ are as follows:
\begin{equation}\label{eqn:ind_rep_branching_rules}
\operatorname{Ind}^{G_{m,n+1}}_{G_{m,n} \times G_{m,1}}S^{\lambda^{[m]}} \otimes S^{\tau^{[m]}} \cong \bigoplus_{\substack{\lambda^{[m]} \subseteq \mu^{[m]} \\ \tau^{[m]} \subseteq \mu^{[m]} }}S^{\mu^{[m]}} 
\end{equation}
\begin{equation}\label{eqn:res_rep_branching_rules}
\operatorname{Res}^{G_{m,n + 1}}_{G_{m,n} \times G_{m, 1}}S^{\mu^{[m]}} \cong \bigoplus_{\substack{\lambda^{[m]} \subseteq \mu^{[m]}\\\tau^{[m]} \subseteq \mu^{[m]}}}S^{\lambda^{[m]}} \otimes S^{\tau^{[m]}}
\end{equation}
\end{theorem}
\begin{proof}
By Theorem \ref{thm:basis_gensymgroup_specht}, we know that $S^{\lambda^{[m]}} \cong \lambda$ where $\lambda$ is an irreducible representation of $G_{m,n}$. Similarly, $\tau$ is an irreducible representations of $G_{m, 1}$. Hence, by Definition \ref{def:outer_product}, $\operatorname{Ind}^{G_{m,n+1}}_{G_{m,n} \times G_{m,1}}S^{\lambda^{[m]}} \otimes S^{\tau^{[m]}}$ is $\lambda \# \tau$. By Theorem \ref{thm:littwoodrichardson}, we can see that the construction for all diagrams of the irreducible components of $\lambda \# \tau$ should contain both $\lambda^{[m]}$ and $\tau^{[m]}$. Hence, $\operatorname{Ind}^{G_{m,n+1}}_{G_{m,n} \times G_{m,1}}S^{\lambda^{[m]}} \otimes S^{\tau^{[m]}}$ would just be the direct sum of the Specht modules that correspond to $\mu^{[m]}$ such that $\lambda^{[m]} \subseteq \mu^{[m]}$ and $\tau^{[m]} \subseteq \mu^{[m]}$ i.e. 
\[\operatorname{Ind}^{G_{m,n+1}}_{G_{m,n}\times G_{m,1}}S^{\lambda^{[m]}} \otimes S^{\tau^{[m]}} \cong \bigoplus_{\substack{\lambda^{[m]} \subseteq \mu^{[m]} \\ \tau^{[m]} \subseteq \mu^{[m]} }}m_\mu S^{\mu^{[m]}}.\]
The multiplicity $m_\mu$ for each summand in the direct sum of Specht modules has to be 1 because the construction of the irreducible components detailed in Theorem \ref{thm:littwoodrichardson} only allows for each irreducible component to be constructed once. Hence, we have shown Equation \ref{eqn:ind_rep_branching_rules}. Equation \ref{eqn:res_rep_branching_rules} can be proven by a standard application of Frobenius reciprocity.
\end{proof}

\section{Future Work}\label{sec:future-work}
In \cite{M20}, Matheau-Raven conjecturized that the eigenvalues of the random transposition shuffle on $G_{m,n}$ may be described by the techniques of lifting vectors.
To prove this conjecture, one has to find correct adding operator and switching operators for the random transposition shuffling on $G_{m,n}$.

In this section, we propose and discuss possible ways of finding such operators. As argued in Matheau-Raven's dissertation, the tabloid notation is inconvenient, and hence it is better to transform to another notation with one-to-one correspondence, \emph{viz}, words.

We here define words for the generalized symmetric group mimicking Definition 2.2.27 in \cite{M20}.
\begin{definition}[Words for The Generalized Symmetric Group]
Define the set 
$$[\overline{n}] = \left(\bigcup_{k \in \{1, \dots, m-1\}} [n^{(k)}]\right) \cup \left(\bigcup_{i \in \{1, \dots, m\}} [n^{\xi_i}]\right),$$ 
where
\[[n^{(k)}] = \{1^{(k)}, \dots, n^{(k)}\}\]
\[[n^{\xi_i}] = \{1^{\xi_i}, \dots, n^{\xi_i}\}\]
We call $[n^{(k)}]$ the unsigned and $[n^{\xi_i}]$ the signed.
We call the set $W^{\overline{n}}$ the set of words lengths $n$ with letters from $[\overline{n}]$, i.e.,
\[w = w_1 w_2 w_3\dots w_n,\]
where
\[w \in W^{\overline{n}} \text{ and } w_i \in [\overline{n}].\]

For the empty word, a string with no letter in there, we denote it as $\omega$.
\end{definition}

\begin{definition}[Map from Partitions to Words]
Let $\lambda^{[m]}$ be a $m-$partition of $n$.
Define a map $w$ that takes a partition and sent to $W^{\overline{n}}$.
Construct it similarly to the function $w$ in \cite{M20} Definition 2.2.27, except for $\lambda^m$, the words are superscripted with their corresponding $\xi_i$, and for $\lambda^1$ to $\lambda^{m-1}$, the words are superscripted with $k$, which represents the partition they are in.

Define the $\operatorname{eval}(w): W^{\overline{n}} \to \lambda^{[m]}$ similarly to the evaluation function in \cite{M20}.
\end{definition}

\begin{lemma}[Word Equivalence]
The word notation satisfies the row equivalence relationship of tableau.
\end{lemma}

\begin{proof}
For $\lambda^1$ to $\lambda^{m-1}$, the superscript of letters record the partition and the integers record row number of each element, and since the row number and partition position are invariant under row equivalence for $\lambda^1$ to $\lambda^{m-1}$, the letters stay the same.

For $\lambda^m$, the superscripts of letters record the $\xi_i$ script and integers row number of each element, and since they are invariant under row equivalence, so the letters stay the same.

Since all the letter stay the same, we know that the word is in one-to-one correspondence  with the tabloids.
\end{proof}

\begin{definition}
The random transposition shuffle on a generalized symmetric group $G_{m,n}$ is driven by the probability distribution $\overline{RT}_{m,n}$
\[\overline{RT}_{m,n}(\sigma) = \begin{cases} 1/mn & \text{if } \sigma = e, \\
1/mn^2 & \text{if } \sigma = \xi_i^k \text{ for } i \in [n], k \in [m-1]\\ 
2/mn^2 & \text{if } \sigma = (ij) \text{ for } i,j \in [n] \text{ with } i < j,\\ 
2/mn^2 & \text{if } \sigma = \xi_i^k \xi_j^k (ij) \text{ for } i,j \in [n] \text{ with } i < j; k \in [m-1]\\ 
0 & \text{otherwise}.\end{cases}\]
\end{definition}

Starting from here, the author provides insights to prove Conjecture 4.5.5 in \cite{M20}.
First, we restate the conjecture as follows:
\begin{conjecture}[{\cite[Conjecture 4.5.5]{M20}}]\label{conj:main}
The eigenvalues for the random transposition shuffle on the generalised symmetric group $G_{m,n}$ are labelled by $m-$partitions of $n$, and may be described by the technique of lifting eigenvectors.
\[\operatorname{eig}(\overline{\lambda}) = \frac{1}{m n^2} (m\abs{\lambda^1} + 2m \operatorname{Diag}(\lambda^1) + \cdots + 2m \operatorname{Diag}(\lambda^m))\]
\end{conjecture}

To prove this theorem, the first step is to write the random transposition shuffle as an element of the vector space $\mathfrak{G}_{m,n} = \mathbb{C}[G_{m.n}]$.

\begin{definition}[$\overline{ART}_{m,n}$]
The random transposition shuffle on $G_{m,n}$ maybe viewed as the following element in $\mathfrak{G}_{m,n}$.
\begin{align}
    \overline{ART}_{m,n} &= m n^2 \sum_{\sigma \in G_{m.n}} \overline{RT}_{m,n}(\sigma) \sigma\\
    &= n \cdot e + \sum_{k \in [m-1]} \sum_{1 \leq i \leq n} \xi_i^k + 2 \sum_{1 \leq i < j \leq n} (ij) + 2 \sum_{k \in [m-1]} \sum_{1 \leq i < j \leq n} \xi_i^k \xi_j^k (ij)
\end{align}
Notice that here, we normalized it by a factor of $m n^2$.
\end{definition}

To justify for our intuition on the lifting property, we calculate the matrix $\overline{ART}_{m,n+1} - \overline{ART}_{m,n}$:
\begin{multline}
    \overline{ART}_{m,n+1} - \overline{ART}_{m,n} \\
    = e + \sum_{k \in [m-1]} \xi_{n+1}^k + 2 \sum_{1 \leq i \leq n} \mqty(i & n+1) + 2 \sum_{k \in [m-1]} \sum_{1 \leq i \leq n} \xi_i^k \xi_{n+1}^k \mqty(i & n+1).
\end{multline}
We find that the resulting terms are only dependent on transpositions involving $(n+1)$.
Hence, intuitively, the lifting eigenvector technique should not only work for hyperoctahedral groups, but also the generalized symmetric groups.

Mimicking the proof of Theorem 4.3.2 in \cite{M20}, which is a hyperoctahedral version of Conjecture \ref{conj:main}, 
one has to be able to construct eigenvectors of $\overline{ART}_{m,n+1}$ when given $\overline{ART}_{m,n}$.
Then, one may restrict the domain to Specht Modules $S^{\overline{\lambda}}$, and compute the change in eigenvalues when constructing new eigenvectors.
The change in eigenvalues may then be computed using the Branching rules introduced in Section \ref{sec:branchingrulesforgensymmgrp}.
One should find that the change in eigenvalues when adding a box to $\lambda^k$ for $k = 1, m-1$ is always the diagonal value of $\lambda^k$, with normalizations.

The steps to a final proof may seems clear, yet the very first step is more complicated than we expected.
To construct the eigenvectors of $\overline{ART}_{m,n+1}$ when given $\overline{ART}_{m,n}$ requires the following conjecture.

\begin{conjecture}\label{conj:operator_group_commutator}
Let $n \in \N$ and $\overline{\lambda} \vdash n$.
For words in $M^{\overline{\lambda}}$ we have the following equalities:
\begin{equation}\label{eqn:signed_art_switching_op_commutator}\overline{ART}_{m,n+1} \circ \Phi_a^m - \Phi_a^m \circ \overline{ART}_{m,n} = 2 \sum_{1 \leq b \leq n} \Phi_b^m \circ \Theta_{b,a}^m,\end{equation}
\begin{equation}\label{eqn:unsigned_art_switching_op_commutator}\overline{ART}_{m,n+1} \circ \Phi_a^k - \Phi_a^k \circ \overline{ART}_{m,n} = 2 \Phi_a^k + 2 \sum_{1 \leq b \leq n} 2 \cdot \Phi_b^k \circ \Theta_{b,a}^k + (\sum_{i = 1}^m \Phi_a^{\xi^i}) \circ \Theta_{b^{\xi^i},a}\end{equation}
where $k$ is from $1$ to $m-1$.
\end{conjecture}

The $\Phi$ and $\Theta$ functions are what we refer as adding and switching operators. Finding a correct version of such operator is crucial to proving the conjecture above. These operators are useful in finding the upper and lower bounds of the eigenvalues of $\lambda^k$. Moreover, we focus on the requirements of the operators to prove Conjecture \ref{conj:main}.

There are two types of operators to be determined, each of which also requires different basic operators to define. The first type is the \emph{adding operator}. Intuitively, suppose $\lambda \vdash n$, then adding operator adds a box to the end of the first row with the number $n+1$.
Adding operator essentially send an element from $M^{\lambda}$ to $M^{\lambda + e_1}$.
Note that the intuition is for the symmetric group $S_n$, so to think about the generalized symmetric group $G_{m,n}$, one has to be careful about the signs of the elements.
Hence, for $G_{m,n}$, one has to define $m$ numbers of adding operators, denoted as $\Phi_{a}^k(w)$ where $w$ is a word in $M^{\overline{\lambda}}$.
This operator means that we add a box at the $a$th row to the $k$th partition $\lambda^{(k)}$ and fill with the value $(n+1)$, where $k$ ranges from $1$ to $m-1$.
For $k = m$, finding the adding operator is crucial, as for all the previous partitions, we add an unsigned element, where for the $m$th partition, we are adding a signed element.
The key here is to find the correct linear combination such that the negative transposition in Equation \eqref{eqn:signed_art_switching_op_commutator} is zero.

For a clear picture of $\Phi_a^k(w)$ where $k$ ranges from $1$ to $m-1$, we here present a rough definition of the basic adding operators:
\begin{align}
    \Phi_a^{\xi^i}(w) &= w a^{\xi^i},\\
    \Phi_a^{k}(w) &= w a^{(k)},
\end{align}
where $i$ ranges from $1$ to $m$ and $k$ ranges from $1$ to $m-1$.

The second type is the \emph{switching operator}.
Suppose $w \in M^{\lambda}$ with $\lambda \vdash n$, and $a, b \in [n]$.
Then the switching operator $\Theta_{b,a}^k(w)$ means summing all of the words with the $i$th position of $w$ being $b$ to $a$.
And $k$ ranges from $1$ to $m-1$ as before.
Similar to the adding operator, when $k$ ranges from $1$ to $m-1$, the switching operator is switching unsigned $a$ with unsigned $b$, yet when either $a$ or $b$ is signed, or both are signed, the switching operator is defined differently.
Nevertheless, the key is the same as before, find a linear combination of the signed switching operators such that the negative transformation of Equation \eqref{eqn:signed_art_switching_op_commutator} is zero.

We present a rough definition of the basic switching operators for a better intuition, as follows:
\begin{align}
    \Theta^{k}_{b,a}(w) &= \sum_{\substack{1 \leq i \leq n \\ w_i = b^{(k)}}} w_1 \dots w_{i-1} a^{(k)} w_{i+1} \dots w_n,\label{eqn:unsigned_switching_op}\\
    \Theta_{b^{\xi^h},a^{\xi^k}} &= \sum_{\substack{1 \leq i \leq n \\ w_i = b^{\xi^h}}} w_1 \dots w_{i-1} a^{\xi^k} w_{i+1} \dots w_n,\label{eqn:signed_switching_op}\\
    \Theta^{\xi^k}_{b,a} &= \sum_{i=1}^{m}\Theta_{b^{\xi^i},a^{\xi^{i+(\text{a phase that involves $k$})}}},\label{eqn:signed_switching_op_plus}\\
    \Theta_{b^{\xi},a^{(k)}} &= \sum_{i=1}^m \Theta_{b^{\xi^i},a^{(k)}}\label{eqn:mix_switching_op},
\end{align}
where $k$ in Equation \eqref{eqn:unsigned_switching_op} ranges from $1$ to $m-1$ and they deal with switches between $a$ and $b$ when they are both unsigned;
the letters $k$ and $h$ in Equation \eqref{eqn:signed_switching_op} ranges from $1$ to $m$ and they deal with switches between $a$ and $b$ when both are signed;
Equation \eqref{eqn:signed_switching_op_plus} defines the switching operator when both $a$ and $b$ are signed with a phase between their signs;
the last equation above deal with switches when $b$ is signed and $a$ is unsigned.

Notice that both operators require finding correct linear combinations for the signed elements, and each operator has $m$ degrees of freedom (because there are $m$ numbers of signs), so the problem has $2m$ degrees of freedom, which is almost impossible to solve.
However, we did find that if the coefficients are roots of unities $\xi_i$ where $i$ ranges from $1$ to $m$, also namely the signs, then such linear combination could satisfy the requirements above.
The only problem left is which root of unity corresponds to which signed element: a task that we leave as an open problem.

Even though we could not give a precise definition of the adding operator and switching operator, we present the format of them as follows:
\begin{align}
    \Phi_a^m(w) &= \sum_{i\in [m]} (\text{The correct choice of the root of unity}) \Phi^{\xi^i}_a(w),\\
    \Theta^m_{b,a} &= \sum_{i=1}^m (\text{The correct choice of the root of unity}) \Theta^{\xi^i}_{b,a}.
\end{align}

Once one finds the correct adding and switching operators, one can prove Conjecture \ref{conj:operator_group_commutator} with the following Lemmas.

\begin{lemma}\label{lem:switching_op_homomorphism}
The switching operators $\Theta_{b^{\xi^i}, a^{(k)}}, \Theta^k_{b,a}$ for $k \in \{1, \ldots, m\}$ are $\mathbb{C}[G_{m,n}]$-module homomorphisms.
\end{lemma}

\begin{lemma}\label{lem:operator_commutator}
The adding and switching operators for $k = \{1, \ldots, m+1\}$satisfy the following equalities:
\begin{equation}\label{eqn:operator_commutator_k}\Phi_b^k \circ \Theta_{b,a}^k = \Theta_{b,a}^k \circ \Phi_b^k - \Phi_a^k,\end{equation}

\begin{equation}\label{eqn:operator_commutator_1}\Phi_b^m \circ \Theta_{a,b}^m = \Theta_{a,b}^m \circ \Phi_b^m - m \cdot \Phi_a^m.\end{equation}
\end{lemma}

The lifting eigenvector technique can provide us with the eigenvalues 
 that can then be used to prove Conjectures~\ref{conjecture:biasedOST} and \ref{conjecture:random}, allowing us to determine the mixing time for the biased OST shuffle and the random transposition shuffle on $G_{m,n}$.

\printbibliography

\end{document}